\documentstyle[12pt,amssymb,amsfonts,twoside]{article}

\begin{titlepage}

\title{Nonexistence of invariant distributions supported on the limit
set} 
\author{Ulrich Bunke and Martin Olbrich} 

\end{titlepage}

\begin{document}

\maketitle
\newcommand{\diag}{{\rm diag}}
\newcommand{\proof}{{\it Proof.$\:\:\:\:$}}
 \newcommand{\dist}{{\rm dist}}
\newcommand{\kaaa}{{\frak k}}
\newcommand{\paaa}{{\frak p}}
\newcommand{\vp}{{\varphi}}
\newcommand{\taaa}{{\frak t}}
\newcommand{\haaa}{{\frak h}}
\newcommand{\R}{{\Bbb R}}
\newcommand{\Hh}{{\bf H}}
\newcommand{\Rep}{{\rm Rep}}
\newcommand{\Hb}{{\Bbb H}}
\newcommand{\Q}{{\Bbb Q}}
\newcommand{\str}{{\rm str}}
\newcommand{\Ind}{{\rm ind}}
\newcommand{\triv}{{\rm triv}}
\newcommand{\Z}{{\Bbb Z}}
\newcommand{\bD}{{\bf D}}
\newcommand{\bF}{{\bf F}}
\newcommand{\tX}{{\tt X}}
\newcommand{\Cliff}{{\rm Cliff}}
\newcommand{\tY}{{\tt Y}}
\newcommand{\tZ}{{\tt Z}}
\newcommand{\tV}{{\tt V}}
\newcommand{\tR}{{\tt R}}
\newcommand{\Fam}{{\rm Fam}}
\newcommand{\Cusp}{{\rm Cusp}}
\newcommand{\bT}{{\bf T}}
\newcommand{\bK}{{\bf K}}
\newcommand{\K}{{\Bbb K}}
\newcommand{\tH}{{\tt H}}
\newcommand{\bS}{{\bf S}}
\newcommand{\bB}{{\bf B}}
\newcommand{\tW}{{\tt W}}
\newcommand{\tF}{{\tt F}}
\newcommand{\bA}{{\bf A}}
\newcommand{\bL}{{\bf L}}
 \newcommand{\bom}{{\bf \Omega}}
\newcommand{\bundle}{{bundle}}
\newcommand{\ch}{{\bf ch}}
\newcommand{\ve}{{\varepsilon}}
\newcommand{\C}{{\Bbb C}}
\newcommand{\gen}{{\rm gen}}
\newcommand{\cTop}{{{\cal T}op}}
\newcommand{\bP}{{\bf P}}
\newcommand{\Naaa}{{\bf N}}
\newcommand{\image}{{\rm image}}
\newcommand{\gaaa}{{\frak g}}
\newcommand{\zaaa}{{\frak z}}
\newcommand{\saaa}{{\frak s}}
\newcommand{\laaa}{{\frak l}}
\newcommand{\stimes}{{\times\hspace{-1mm}\bf |}}
\newcommand{\ausg}{{\rm end}}
\newcommand{\bff}{{\bf f}}
\newcommand{\maaa}{{\frak m}}
\newcommand{\aaaa}{{\frak a}}
\newcommand{\naaa}{{\frak n}}
\newcommand{\brr}{{\bf r}}
\newcommand{\res}{{\rm res}}
\newcommand{\Aut}{{\rm Aut}}
\newcommand{\Pol}{{\rm Pol}}
\newcommand{\Tr}{{\rm Tr}}
\newcommand{\cT}{{\cal T}}
\newcommand{\dom}{{\rm dom}}
\newcommand{\db}{{\bar{\partial}}}
\newcommand{\g}{{\gaaa}}
\newcommand{\cZ}{{\cal Z}}
\newcommand{\cH}{{\cal H}}
\newcommand{\cM}{{\cal M}}
\newcommand{\interi}{{\rm int}}
\newcommand{\singsupp}{{\rm singsupp}}
\newcommand{\cE}{{\cal E}}
\newcommand{\ccR}{{\cal R}}
\newcommand{\cV}{{\cal V}}
\newcommand{\cY}{{\cal Y}}
\newcommand{\cW}{{\cal W}}
\newcommand{\cI}{{\cal I}}
\newcommand{\cC}{{\cal C}}
\newcommand{\mod}{{\rm mod}}
\newcommand{\cK}{{\cal K}}
\newcommand{\cA}{{\cal A}}
\newcommand{\cEp}{{{\cal E}^\prime}}
\newcommand{\cU}{{\cal U}}
\newcommand{\Hom}{{\mbox{\rm Hom}}}
\newcommand{\vol}{{\rm vol}}
\newcommand{\cO}{{\rm \bf o}}
\newcommand{\End}{{\mbox{\rm End}}}
\newcommand{\Ext}{{\mbox{\rm Ext}}}
\newcommand{\rk}{{\rm rank}}
\newcommand{\im}{{\mbox{\rm im}}}
\newcommand{\sign}{{\rm sign}}
\newcommand{\spann}{{\mbox{\rm span}}}
\newcommand{\symm}{{\mbox{\rm symm}}}
\newcommand{\cF}{{\cal F}}
\newcommand{\cD}{{\cal D}}
\newcommand{\Ree}{{\rm Re }}
\newcommand{\Res}{{\mbox{\rm Res}}}
\newcommand{\Imm}{{\rm Im}}
\newcommand{\inter}{{\rm int}}
\newcommand{\clo}{{\rm clo}}
\newcommand{\tg}{{\rm tg}}
\newcommand{\ee}{{\rm e}}
\newcommand{\Li}{{\rm Li}}
\newcommand{\cN}{{\cal N}}
 \newcommand{\conv}{{\rm conv}}
\newcommand{\op}{{\mbox{\rm Op}}}
\newcommand{\tr}{{\mbox{\rm tr}}}
\newcommand{\cs}{{c_\sigma}}
\newcommand{\ctg}{{\rm ctg}}
\newcommand{\degg}{{\mbox{\rm deg}}}
\newcommand{\Ad}{{\mbox{\rm Ad}}}
\newcommand{\ad}{{\mbox{\rm ad}}}
\newcommand{\codim}{{\rm codim}}
\newcommand{\Gr}{{\mathrm{Gr}}}
\newcommand{\coker}{{\rm coker}}
\newcommand{\id}{{\mbox{\rm id}}}
\newcommand{\ord}{{\rm ord}}
\newcommand{\nat}{{\Bbb  N}}
\newcommand{\supp}{{\rm supp}}
\newcommand{\sing}{{\mbox{\rm sing}}}
\newcommand{\spec}{{\mbox{\rm spec}}}
\newcommand{\Ann}{{\mbox{\rm Ann}}}
\newcommand{\aca}{{\aaaa_\C^\ast}}
\newcommand{\acag}{{\aaaa_{\C,good}^\ast}}
\newcommand{\acage}{{\aaaa_{\C,good}^{\ast,extended}}}
\newcommand{\tck}{{\tilde{\ck}}}
\newcommand{\tnk}{{\tilde{\ck}_0}}
\newcommand{\ceep}{{{\cal E}(E)^\prime}}
 \newcommand{\ncE}{{{}^\naaa\cE}}
 \newcommand{\Or}{{\rm Or }}
\newcommand{\Diff}{{\cal D}iff}
\newcommand{\cB}{{\cal B}}
\newcommand{\hc}{{{\cal HC}(\gaaa,K)}}
\newcommand{\hcma}{{{\cal HC}(\maaa_P\oplus\aaaa_P,K_P)}}
\def\imath{{\rm i}}
\newcommand{\vsl}{{V_{\sigma_\lambda}}}
\newcommand{\czg}{{\cZ(\gaaa)}}
\newcommand{\csl}{{\chi_{\sigma,\lambda}}}
\newcommand{\cR}{{R}}
\def\hB{\hspace*{\fill}$\Box$ \newline\noindent}
\newcommand{\varho}{\varrho}
\newcommand{\ind}{{\rm index}}
\newcommand{\Indu}{{\rm Ind}}
\newcommand{\Fin}{{\mbox{\rm Fin}}}
\newcommand{\cS}{{S}}
\newcommand{\orig}{{\cal O}}
\def\hB{\hspace*{\fill}$\Box$ \\[0.5cm]\noindent}
\newcommand{\cL}{{\cal L}}
 \newcommand{\cG}{{\cal G}}
\newcommand{\npci}{{\naaa_P\hspace{-1.5mm}-\hspace{-2mm}\mbox{\rm coinv}}}
\newcommand{\pki}{{(\paaa,K_P)\hspace{-1.5mm}-\hspace{-2mm}\mbox{\rm inv}}}
\newcommand{\mki}{{(\maaa_P\oplus \aaaa_P, K_P)\hspace{-1.5mm}-\hspace{-2mm}\mbox{\rm inv}}}
\newcommand{\Mat}{{\rm Mat}}
\newcommand{\npi}{{\naaa_P\hspace{-1.5mm}-\hspace{-2mm}\mbox{\rm inv}}}
\newcommand{\ngp}{{N_\Gamma(\pi)}}
\newcommand{\gbg}{{\Gamma\backslash G}}
\newcommand{\gkm}{{ Mod(\gaaa,K) }}
\newcommand{\ggkm}{{  (\gaaa,K) }}
\newcommand{\pkm}{{ Mod(\paaa,K_P)}}
\newcommand{\ppkm}{{  (\paaa,K_P)}}
\newcommand{\makm}{{Mod(\maaa_P\oplus\aaaa_P,K_P)}}
\newcommand{\mmakm}{{ (\maaa_P\oplus\aaaa_P,K_P)}}
\newcommand{\cP}{{\cal P}}
\newcommand{\gm}{{Mod(G)}}
\newcommand{\gk}{{\Gamma_K}}
\newcommand{\La}{{\cal L}}
\newcommand{\cug}{{\cU(\gaaa)}}
\newcommand{\cuk}{{\cU(\kaaa)}}
\newcommand{\dc}{{C^{-\infty}_c(G) }}
\newcommand{\gdk}{{\gaaa/\kaaa}}
\newcommand{\dgkm}{{ D^+(\gaaa,K)-\mbox{\rm mod}}}
\newcommand{\dgm}{{D^+G-\mbox{\rm mod}}}
\newcommand{\vect}{{\C-\mbox{\rm vect}}}
 \newcommand{\cig}{{C^{ \infty}(G)_{K} }}
\newcommand{\gami}{{\Gamma\hspace{-1.5mm}-\hspace{-2mm}\mbox{\rm inv}}}
\newcommand{\cQ}{{\cal Q}}
\newcommand{\mmap}{{Mod(M_PA_P)}}
\newcommand{\bbbz}{{\bf Z}}
 \newcommand{\cX}{{\cal X}}
\newcommand{\bH}{{\bf H}}
\newcommand{\pr}{{\rm pr}}
\newcommand{\bX}{{\bf X}}
\newcommand{\bY}{{\bf Y}}
\newcommand{\bZ}{{\bf Z}}
\newcommand{\bV}{{\bf V}}

\newtheorem{prop}{Proposition}[section]
\newtheorem{lem}[prop]{Lemma}
\newtheorem{ddd}[prop]{Definition}
\newtheorem{theorem}[prop]{Theorem}
\newtheorem{kor}[prop]{Corollary}
\newtheorem{ass}[prop]{Assumption}
\newtheorem{con}[prop]{Conjecture}
\newtheorem{prob}[prop]{Problem}
\newtheorem{fact}[prop]{Fact}

\tableofcontents
\parskip3ex
\section{Statement of the result}
Let $X$ be a symmetric space of negative curvature. 
Then $X$ either belongs to one of the three families of real, complex, or
quaternionic hyperbolic spaces, or it is the Cayley hyperbolic plane.
 
Let $G$ be a connected linear Lie group which finitely covers the isometry
group of $X$. Furthermore, let $\Gamma\subset G$ be a discrete subgroup. We
assume that $\Gamma$ is geometrically finite. We refer to Definition \ref{t67}
for a precise explanation of this notion. If $X$ is a real  hyperbolic space,
then $\Gamma$ is geometrically finite iff it admits a fundamental domain with
finitely many totally geodesic faces. In the other cases the definition is
more complicated. Essentially, $\Gamma$ is geometrically finite if the
corresponding locally symmetric space $\Gamma\backslash X$ has finitely many
cusps and can be compactified by adding a geodesic boundary and closing the
cusps. In particular, if $\Gamma$ is cocompact, or  convex-cocompact, or the
locally symmetric space $\Gamma\backslash X$ has finite volume, then
$\Gamma$ is geometrically finite.  

We adjoin the geodesic boundary $\partial X$ to $X$ and obtain a 
compact manifold with boundary $\bar X:=X\cup\partial X$. 
A point of $\partial X$ is an equivalence class of geodesic rays where two
rays are in the same class if they run in bounded distance to each other.
The action of $G$ extends naturally to $\bar X$. Let $\Lambda_\Gamma\subset
\partial X$ denote the limit set of $\Gamma$. It is defined as the set of
accumulation points of any orbit $\Gamma\cO$ in $\bar X$ for $\cO\subset X$.

We consider a $G$-equivariant irreducible complex vector bundle $V\rightarrow
\partial X$ and a finite-dimensional representation $(\vp,V_\vp)$ of $\Gamma$.
Furthermore, by $\Lambda\rightarrow \partial X$
we denote the $G$-equivariant bundle of densities on $\partial X$.
To $V$ we associate the $G$-equivariant bundle  
$$\tilde V:=\Hom(V,\Lambda)\ .$$  
The space $C^{-\infty}(\partial X,V)$ of
distribution sections of $V$ is then, by definition, the topological dual of
$C^{\infty}(\partial X, \tilde V)$. 
We define the space of invariant distribution sections of $V$ with
twist $\vp$ by \begin{ddd}
$$I(\Gamma,V,\vp):=\left( C^\infty(\partial X,V)\otimes V_\vp \right)^\Gamma \
.$$
\end{ddd} 

Next we introduce some real quantities which represent growth properties of the
geometric objects introduced so far.
We first define the number $\rho\in\R$
which is  a measure of the volume growth of the symmetric space $X$. We use
this number in order to normalize the critical exponents below. 
Let $\cO$ be
any point of $X$, and let $B(r,\cO)$ denote the ball of radius $r$ centered at
$\cO$. 
\begin{ddd} $$\rho:=\frac12 \lim_{r\to\infty} \frac{\log\:\vol
\:B(r,\cO)}{r}\ .$$ 
\end{ddd}  
 The growth of the action of $G$ on the
bundle $V$ is measured by the quantity $s(V)\in\R$. Note that $\Lambda$ is the
complexification of a real orientable line bundle. It is therefore trivial if
considered merely as a vector bundle, but it is not trivial as a
$G$-equivariant bundle. The bundle $\Lambda$ can be represented by a cocycle
of positive transition functions. If $\alpha\in\C$, then raising the
transition functions to the power $\alpha$, we obtain a new cocycle which
represents the $G$-equivariant bundle $\Lambda^\alpha$.

\begin{ddd}
$s(V)$is  defined as the unique number such that $V\otimes
\Lambda^{s(V)}\cong V^\sharp$ as $G$-equivariant bundles, where
$V^\sharp$ denotes the complex conjugate bundle of $\tilde V$.  
\end{ddd}
For example, if $V=\partial X\times\C$ is the trivial bundle, then $s(V)=1$.
 More generally, $s(\Lambda^\alpha)=1-2\Ree(\alpha)$. 
 
The normalized growth of $\Gamma$ is expressed by the critical exponent
\begin{ddd}
$$d_\Gamma:=\frac{1}{\rho} \inf\{\nu\:|\: \sum_{g\in\Gamma}
\dist(g\cO,\cO)^{-\rho-\nu}<\infty\}\ .$$
\end{ddd}
This definition is independent of the choice of $\cO\in X$.
Since $\Gamma$ is discrete and infinite we have $d_\Gamma\in (-1,1]$.

The exponent $d_\vp$ is a measure for the growth of $\vp$. It is
defined by
\begin{ddd}
$$d_\vp:=\frac{1}{\rho} \inf\{\nu\:|\:
\sup_{g\in\Gamma}\|\vp(g)\|\dist(g\cO,\cO)^{-\nu} < \infty \}\ ,$$
\end{ddd}
where we have fixed any norm $\|.\|$ on $\End(V_\vp)$ and any point $\cO\in X$.
Since $\Gamma$ is finitely generated, we have $d_\vp<\infty$.

A cusp of $\Gamma$ is, by definition, a $\Gamma$-conjugacy class $[P]_\Gamma$
of proper parabolic subgroups $P\subset G$ such that $\Gamma\cap P$ is infinite
and $\pi(\Gamma\cap P)\subset L$ is precompact, where $\pi$ is the projection
onto the semisimple quotient $L$ given by the sequence  $$0\rightarrow
N\rightarrow P\stackrel{\pi}{\rightarrow} L\rightarrow 0$$ with $N\subset P$
denoting the unipotent radical of $P$. Note that if $[P]_\Gamma$ is a cusp of
$\Gamma$, then $\Gamma_P:=\Gamma\cap P$ again satisfies our assumptions.
The limit set of $\Gamma_P$ consists of the unique fixed point
$\infty_P\subset \partial X$ of $P$. Since $\Gamma_P$ acts properly on
$\Omega_{\Gamma_P}:=\partial X\setminus \{\infty_P\}$ and
$\Gamma_P\backslash (\Lambda_\Gamma\setminus \{\infty_P\})\subset
\Gamma_P\backslash \Omega_{\Gamma_P}$ is compact (see Lemma \ref{hilfe}) we can
choose a smooth function $\chi^{\Gamma_P}$ on $\Omega_{\Gamma_P}$ such that
$\supp(\chi^{\Gamma_P})\cap \Lambda_\Gamma$ is a compact subset of
$\Omega_{\Gamma_P}$, $\{\supp(g^*\chi^{\Gamma_P})\}_{g\in\Gamma_P}$ is a
locally finite covering of $\Omega_{\Gamma_P}$, and $\sum_{g\in\Gamma_P}g^*
\chi^{\Gamma_P}\equiv 1$.   Assume that $s(V)>d_\Gamma+d_\vp$. 
\begin{ddd}\label{str}
We say that $f\in I(\Gamma,V,\vp)$ is strongly supported on the limit set if
\begin{enumerate}
\item $f$ is supported on the limit set as a distribution.
\item For any $h\in V^\infty(\partial X,\tilde V)\otimes\tilde V_\vp$ and cusp
$[P]_\Gamma$ of $\Gamma$ we have 
$$\langle f,h\rangle =\sum_{g\in \Gamma_P} \langle \chi^{\Gamma_P}
f_{|\Omega_{\Gamma_P}}, \tilde \vp(g)^{-1}g^* h\rangle \ .$$
\end{enumerate}
\end{ddd}
In order to see that the second condition is well-defined note that
$\supp(\chi^{\Gamma_P} f_{|\Omega_{\Gamma_P}})\subset
\supp(\chi^{\Gamma_P})\cap \Lambda_\Gamma$ is a compact subset of
$\Omega_{\Gamma_P}$. Therefore the pairing $\langle \chi^{\Gamma_P}
f_{|\Omega_{\Gamma_P}}, \tilde \vp(g)^{-1}g^* h\rangle $ is defined.
The sum converges because of our assumption
$s(V)>d_\Gamma+d_\vp\ge d_{\Gamma_P}+d_\vp$, which implies that
$\sum_{g\in \Gamma_P}  \tilde \vp(g)^{-1}g^* h_{|\Omega_{\Gamma_P}}$
converges in the space of smooth functions. In fact, the argument proving
\cite{bunkeolbrich982}, Lemma 4.2, applies in the more general case when
$\Gamma$ is merely geometrically finite.
In Lemma \ref{independent} 
we will verify that this definition is independent of the choice of
$\chi^{\Gamma_P}$.

In \cite{bunkeolbrich011} and
\cite{bunkeolbrich02} we have expressed the condition "strongly supported on
the limit set" in the form $res^\Gamma(f)=0$. While this definition works for
all values of $s(V)$ there we must assume that $f$ is "deformable". Because in
the present paper we are in the "domain of convergence" we can use the simpler
and more general definition above. 

\begin{ddd}
By 
$I_{\Lambda_\Gamma}(\Gamma,V,\vp)$ we denote the subspace of all $f\in I(\Gamma,V,\vp)$
which are strongly supported on the limit set.
\end{ddd}

The main result of the present paper can now be formulated as follows.
\begin{theorem}\label{main}
If $s(V)>d_\Gamma+d_\vp+\max_{[P]_\Gamma
}(0,d_{\Gamma_P}-d_\Gamma+1)$ (where the maximum is taken over all
cusps of $\Gamma$), then
$I_{\Lambda_\Gamma}(\Gamma,V,\vp)=0$. \end{theorem}

Let us note the following special case which was already shown in
\cite{bunkeolbrich982}, Thm 4.7. The group $\Gamma$ is called convex cocompact
if it acts freely and cocompactly on $\bar X\setminus \Lambda_\Gamma$. In this
case $\Gamma$ has no cusps and
$I_{\Lambda_\Gamma}(\Gamma,V,\vp)$ is just the space invariant distribution
sections of $V$ with twist $\vp$ which are supported on $\Lambda_\Gamma$.
\begin{kor}If $\Gamma$ is convex cocompact
and $s(V)>d_\Gamma+d_\vp$, then $I_{\Lambda_\Gamma}(\Gamma,V,\vp)=0$.
\end{kor}

Back to the general case of a geometrically finite discrete group  let 
$1$ be the trivial representation of $\Gamma$. Then we have
$d_1=0$. In the place of $V$ we consider $\Lambda^{\frac{1-d_\Gamma}{2}}$. 
Note that $s(\Lambda^{\frac{1-d_\Gamma}{2}})=d_\Gamma=d_\Gamma+d_\vp$.
The space $I_{\Lambda_\Gamma}(\Gamma,\Lambda^{\frac{1-d_\Gamma}{2}},1)$ is
spanned by the Patterson-Sullivan measure \cite{patterson76},
\cite{sullivan79}, \cite{corlette90},  \cite{corletteiozzi99}, hence $\dim 
I_{\Lambda_\Gamma}(\Gamma,\Lambda^{\frac{1-d_\Gamma}{2}},1)=1$. 
Here we must use the definition of the condition "strongly supported on the
limit set" in terms of $res^\Gamma$ given in \cite{bunkeolbrich011}. 

In order to
construct some twisted examples we consider a finite-dimensional $M$-spherical
representation $(\pi,V_\pi)$ of $G$. Here $(\pi,V_\pi)$ is called
$M$-spherical, if for any parabolic subgroup $P\subset G$ there exists a
vector $0\not=v\in V_\pi$ and a character $\chi:P\rightarrow \R$ such that
$\pi(p)v=\chi(p) v$ for all $p\in P$. There is a natural inclusion
$$I_{\Lambda_\Gamma}(\Gamma,\Lambda^{\frac{1-d_\Gamma}{2}},1)\hookrightarrow
I_{\Lambda_\Gamma}(\Gamma,\Lambda^{\frac{1-d_\Gamma-d_\pi}{2}},\pi)$$
showing that
$I_{\Lambda_\Gamma}(\Gamma,\Lambda^{\frac{1-d_\Gamma-d_\pi}{2}},\pi)\not=0$.
On the other hand,
$s(\Lambda^{\frac{1-d_\Gamma-d_\pi}{2}})=d_\Gamma+d_\pi$.

These examples show that our estimate can not be improved in general for
convex cocompact $\Gamma$. On the other hand, even for geometrically finite
$\Gamma$ we do not know any counterexample to the assertion that
already $s(V)>d_\Gamma+d_\vp$ implies that $I_{\Lambda_\Gamma}(\Gamma,V,\vp)=0$.

\section{Geometry of geometrically finite discrete subgroups}

If $\Gamma\subset G$ is a discrete subgroup and $\Lambda_\Gamma$ is its limit
set, then $\Gamma$ acts on $\bar X\setminus\Lambda_\Gamma$ properly
discontinuously. 
Let $\bar Y_\Gamma$ denote the manifold
with boundary $\bar Y_\Gamma:=\Gamma\backslash( \bar X\setminus 
\Lambda_\Gamma)$. 
 If $[P]_\Gamma$ is a cusp of $\Gamma$, then   
we form the manifold with boundary $\bar Y_{\Gamma_P}:= \Gamma_P\backslash 
(\bar X\setminus \{\infty_P\})$.

\begin{ddd}\label{t67}
 The group $\Gamma$ is called geometrically finite if the following conditions hold:
\begin{enumerate}
\item $\Gamma$ has finitely many cusps.
\item There is a bijection between the set of ends of $\bar Y_\Gamma$ and 
and the set of cusps of $\Gamma$. 
\item If $[P]_\Gamma$ is a cusp of $\Gamma$, then there exists a
representative $\bar Y_P$ of the corresponding  end of $\bar Y_\Gamma$
and  embedding $e_P: \bar Y_P  \rightarrow \bar Y_{\Gamma_P}$
which is isometric in the interior such
that its image $e_P(\bar Y_P)$ represents the end of $\bar Y_{\Gamma_P}$.
\end{enumerate}
\end{ddd}

\begin{lem}\label{hilfe}
If $[P]_\Gamma$ is a cusp of $\Gamma$, then $\Gamma_P\backslash
(\Lambda_\Gamma\setminus\{\infty_P\})$ is a compact subset of
$\Gamma_P\backslash \Omega_{\Gamma_P}$.\end{lem}
\proof
It suffices to show that $\Gamma_P\backslash
(\Lambda_\Gamma\setminus\{\infty_P\})$ is compact in $\bar Y_{\Gamma_P}$.
Note that $(\Lambda_\Gamma\setminus\{\infty_P\})$ is closed in $\bar
X\setminus \{\infty_P\}$. Therefore,
$\Gamma_P\backslash
(\Lambda_\Gamma\setminus\{\infty_P\})$ is closed in $\bar Y_{\Gamma_P}$.
Furthermore, it is contained in the compact set $\bar Y_{\Gamma_P}\setminus
e_P(\bar Y_P)$ (note that $\bar Y_P$ is open). The assertion now follows. \hB

Let $\cO\in X$ be any point. We consider the Dirichlet domain 
$F\subset X$ of $\Gamma$ with respect to $\cO$. It is a fundamental domain
given by  $$F:=\{x\in X|\dist(x,\cO)\le \dist(hx,\cO)\:\forall h\in\Gamma\}\
.$$  If $[P]_\Gamma$ is a cusp of $\Gamma$, then let $\chi^{\Gamma_P}$ be the
cut-off function introduced before  Definition \ref{str}. 

\begin{lem} \label{dec} 
We can decompose $F$ as $F_0\cup F_1\cup\dots
F_r$, where $r$ is the number of cusps $[P_i]_\Gamma$, $i=1,\dots r$, of
$\Gamma$, and the subsets $F_i$ satisfy
\begin{enumerate}
\item The closure of $F_0$ in $\bar
X\setminus \Lambda_\Gamma$ is compact.
\item  $\overline{\Gamma_{P_i}  F_i}\cap
(\Lambda_\Gamma\setminus \{\infty_{P_i}\})\cap
\supp(\chi^{\Gamma_{P_i}})=\emptyset$ for $i=1,\dots,r$. \end{enumerate}
\end{lem}
\proof
By $\bar Y_0$ we denote the compact subset 
$\bar Y_\Gamma\setminus \bigcup_{i=1,\dots r} \bar Y_{P_i}$
of $\bar Y_\Gamma$. Then we define $F_0:=\Gamma \bar Y_0 \cap F$,
where $\Gamma \bar Y_0$ denotes the preimage of $\bar Y_0$ under the
projection $\bar X\setminus\Lambda_\Gamma\rightarrow \bar Y_\Gamma$.
By definition, $\bar F_0 \subset (\bar X\setminus
\Lambda_\Gamma)$. 

For $i=1,\dots r$ we define $F_i:=\Gamma \bar Y_{P_i}\cap F$. 
We then have $\overline{\Gamma_{P_i}  F_i}\cap
(\Lambda_\Gamma\setminus \{\infty_{P_i}\})\cap
\supp(\chi^{\Gamma_{P_i}})=\emptyset$ since the contrary this would imply
$e(\bar Y_{P_i}) \cap \Gamma_{P_i}\backslash
(\Lambda_\Gamma\setminus\{\infty_{P_i}\})\not=\emptyset$.
\hB

\begin{lem}\label{independent}
Let $\chi_1,\chi_2$ be two choices of the cut-off function $\chi^{\Gamma_P}$
in Definition \ref{str}.
Then $$\sum_{g\in \Gamma_P} \langle \chi_1
f_{|\Omega_{\Gamma_P}}, \tilde \vp(g)^{-1}g^* h\rangle=\sum_{g\in \Gamma_P}
\langle \chi_2 f_{|\Omega_{\Gamma_P}}, \tilde \vp(g)^{-1}g^* h\rangle\ , $$
where $[P]_\Gamma$, $f$, and $h$ are as \ref{str}.
\end{lem}
\proof
The estimates given in the proof of
\cite{bunkeolbrich982}, Lemma 4.2, show that all sums below converge
absolutely. This justifies the resummations in the following computation. In
the first and the last equality we use the $\Gamma$-invariance of $f$.
\begin{eqnarray*} \sum_{g\in \Gamma_P} \langle \chi_1
f_{|\Omega_{\Gamma_P}}, \tilde \vp(g)^{-1}g^* h\rangle&=&\sum_{g\in \Gamma_P}
\langle g^*\chi_1 f_{|\Omega_{\Gamma_P}},   h\rangle\\
&=&\sum_{g\in \Gamma_P}\sum_{l\in\Gamma_P}
\langle g^*\chi_1 l^*\chi_2 f_{|\Omega_{\Gamma_P}},   h\rangle\\
&=&\sum_{l\in \Gamma_P}\sum_{g\in\Gamma_P}
\langle g^*\chi_1 l^*\chi_2 f_{|\Omega_{\Gamma_P}},   h\rangle\\
&=&\sum_{l\in \Gamma_P}
\langle  l^*\chi_2 f_{|\Omega_{\Gamma_P}},   h\rangle\\
&=&\sum_{g\in \Gamma_P}
\langle \chi_2 f_{|\Omega_{\Gamma_P}}, \tilde \vp(g)^{-1}g^* h\rangle\ .
\end{eqnarray*}
\hB

\section{Proof of Theorem \ref{main}}

We adapt the argument of the proof of \cite{bunkeolbrich982}, Thm.4.7
given there in the special case of a convex-cocompact group $\Gamma$ to
the present situation where $\Gamma$ is geometrically finite.

\begin{ddd}
We call the bundle $V$ spherical, if $V=\Lambda^{\frac{1-t(V)}{2}}$ for some
$t(V)\in \C$.
\end{ddd}
Note that $\Ree\:t(V)=s(V)$.
We first show the following special case. 
\begin{prop}\label{poo}
Theorem \ref{main} is true if $V$ is spherical.
\end{prop}
\proof
Let $f\in I_{\Lambda_\Gamma}(\Gamma,V,\vp)$.
Then we must show that $\langle f,h\rangle=0$ for any $h\in C^\infty(\partial
X,\tilde V)\otimes \tilde V_\vp$. 

\begin{lem}
If $\Gamma$ does not contain any hyperbolic element, then $\langle
f,h\rangle=0$. \end{lem}
\proof
If $\Gamma$ does not contain any hyperbolic element, then $\Gamma=\Gamma_P$ for
the unique cusp $[P]_\Gamma$ of $\Gamma$. Since $f$ is supported on
$\Lambda_\Gamma=\{\infty_P\}$ as a distribution we have
we have $f_{|\Omega_{\Gamma_P}}=0$. This implies
$\langle f,h\rangle =\sum_{g\in \Gamma_P} \langle \chi^{\Gamma_P}
f_{|\Omega_{\Gamma_P}},\tilde\vp(g)^{-1}g^* h\rangle=0$.\hB

It remains to consider the case
that $\Gamma$ contains a hyperbolic element which we will denote by $g_0$. 
 
\begin{lem}\label{vam1}
If $\Gamma$ does  contain a hyperbolic element, say $g_0$, then $\langle
f,h\rangle=0$.
\end{lem}
\proof
Let $b_\pm\in \partial X$ denote the attracting and repelling
fixed points of $g_0$. We can write $h=h_++h_-$ such that
$h_\pm$ vanishes in a neighbourhood of $b_\mp$.
It suffices to show that
$\langle f,h\rangle=0$ for any $h$ which vanishes in a neighbourhood
of say $b_+$.

We fix the origin $\cO\in X$  
such that $\cO$ is on the unique geodesic connecting
$b_-$ with $b_+$.
Let $\tilde F\subset X$ be the Dirichlet domain 
 of $\Gamma$ with respect to this choice of the origin. Furthermore,
let $\bar{\tilde F}$ be the closure of $\tilde F$ in $X\setminus
\Lambda_\Gamma$.
The Dirichlet domain $D_{<g_0>}$
with respect to $\cO$
of the group $<g_0>$  generated by $g_0$ separates 
$X\setminus D_{<g_0>}$ into two connected components
$X_+$ and $X_-$. Let $\partial X_\pm:=\bar X_\pm\cap \partial X$.
We can assume that 
$b_\pm \in \partial X_\pm$.

Replacing $\cO$, if necessary, by $g_0^j \cO$, $j\in \nat_0$ sufficiently
large, we can assume that
$\supp(h)\subset  \partial X_-$.
Then we define $F:=g_0^{i} \tilde F$, $\bar F:= g_0^{i} \bar{\tilde F}$,
where we choose $i\in\nat_0$  sufficiently large such that
$F\subset  X_+$.

We use the polar coordinates $(a,k)\in \R_+\times \partial X$ in
order to parametrize points  $x\in X\setminus \{\cO\}$
such that $a(x)=\exp(\dist(\cO,x))$ and $k(x)\in\partial X$ is represented  by
the geodesic ray through $x$ starting in $\cO$.    
Using these coordinates we extend $h$ to the interior of $X$ setting
$\tilde h(x)=\chi(a(x))h(k(x))$, where $\chi\in C^\infty(R_+)$
is some cut-off function which is equal to one near infinity and vanishes
for $a<1$. Note that $\supp(\tilde h)\subset X\setminus X_+$.
 
Note that by our assumption $V$ is spherical and $s(V)>0$. Therefore, the
Poisson transformation $$P:C^{-\infty}(\partial X,V)\rightarrow C^\infty(X)$$
is injective (we refer to \cite{bunkeolbrich982} and the literature cited
therein (e.g.  \cite{schlichtkrull84}) for a definition of the Poisson
transformation and its properties). We use the same symbol $P$ in order to
denote  the extension of the Poisson transform to the tensor product by
$V_\vp$.

Using the polar coordinates we pull-back the volume form of the unit sphere
in $T_\cO X$ to $\partial X$ and thus obtain a volume form $dk$ on $\partial
X$.  Then the inverse of the Poisson transformation is given by the following
limit formula
$$
\langle f,h\rangle =c_1 \lim_{a\to\infty}
a^{\rho(1-t(V))} \int_{\partial X} \langle P f(a(x),k), h(k)\rangle dk
$$
for some constant $c_1$.  
Using the fact that for large $a$ the volume form $dx$ can be written as
$dx=c_2 a^{2\rho}da dk+O(a^{2\rho-1})$ we deduce 
 $$\langle f,h\rangle = c \lim_{a\to\infty} a^{-\rho(1+t(V))}\int_{\{x\in
G|a\le a(x)\le a_0a\}} \langle Pf(x),\tilde h(x)\rangle dx\ .$$
where $c$ depends on $c_1$, $a_0>1$, and $c_2$.
We now employ the covering of $X$ by translates of the fundamental domain
$gF$, $g\in\Gamma$, and the $\Gamma$-invariance of $Pf$
$$Pf(gx)=\vp(g)Pf(x)\ .$$ 
We get 
$$ \langle f,h\rangle =  c\lim_{a\to\infty} a^{-\rho(1+t(V))}\sum_{g\in
\Gamma} \int_{\{x\in F |a\le a(gx)\le a_0a\}}
\langle \vp(g)P f(x),\tilde h(gx)\rangle dx\ .
$$
Since $\supp(\tilde h)\subset
X\setminus X_+$ 
we have  $gF\cap X\setminus X_+\not=\emptyset$
if $g\in\Gamma$ contributes to the sum above.
The triangle inequality for $X$ gives $a(x)a(g)\ge a(gx)$, where we write
$a(g)$ for $a(g\cO)$.  
We will also need the following converse version of the triangle inequality.
\begin{lem}\label{lemm3}
There exists $a_1\in \R_+$ such that for all $g\in \Gamma$  with
$gF\cap  X\setminus X_+\not=\emptyset$  
and $x\in F$ we have $a(g)a(x)\le a_1 a(gx)$.
\end{lem}
We postpone the proof of this lemma and continue the argument for
Lemma \ref{vam1}. Using \ref{lemm3} we obtain
\begin{eqnarray*}
\{x\in F\:|\:a\le a(gx)\le a_0 a\}  &\subset&
\{x\in F\:|\:a\le a(x)a(g)\le a_1a\} \\
&=&\{x\in F\:|\:aa(g)^{-1}\le a(x)\le a_1aa(g)^{-1}\} 
\end{eqnarray*}
for all $g\in \Gamma$ with $gF\cap  X\setminus X_+\not=\emptyset$.
Taking into account that $\tilde h$ is
bounded and that for given $\epsilon>0$ there exists a constant $C_0$ such that
for all $g\in\Gamma$ we have 
 $\|\vp(g)\| \le C_0 a(g)^{\rho(d_\vp+\epsilon)}$
we obtain \begin{eqnarray}
\lefteqn{
|\int_{\{x\in F|a\le a(gx)\le a_0a\} }
\langle \vp(g)P f(x),\tilde h(gx)\rangle dx|
}\nonumber\\
&\le& C_1  a(g)^{\rho(d_\vp+\epsilon)}
\int_{\{x\in F|aa(g)^{-1}\le a(x)\le a_1aa(g)^{-1}\}} |Pf(x)| dx\
,\label{sumabobe} \end{eqnarray}
where $C$ is independent of $g\in\Gamma$ and $a\in\R_+$.
In order to proceed further we employ the following crucial estimate.
\begin{lem}\label{claimm}
There is a constant $C_1$ such that
$$\int_{\{x\in F|b\le a(x)\le a_1 b\} } |Pf(x)| dx\le C_1 
b^{\rho(s(V)+1-\mu)}$$
for all sufficiently small $\mu>0$ and all $b\ge 1$.
\end{lem}
We again postpone the proof of Lemma \ref{claimm} and continue
with the proof of Lemma \ref{vam1}.
If we insert the estimate claimed in Lemma \ref{claimm} into (\ref{sumabobe})
and sum over $\Gamma$, then we obtain
\begin{eqnarray*}
\lefteqn{
\sum_{g\in\Gamma} |\int_{\{x\in F|a\le a(gx)\le a_0a\} }
\langle \vp(g)P f(x),\tilde h(gx)\rangle dx|
}\\
&\le&
C_2\sum_{g\in\Gamma} a(g)^{\rho(d_\vp+\epsilon+\mu-s(V)-1)}
a^{\rho(s(V)+1-\mu)}\ , \end{eqnarray*}
where $C_3$ is independent of $a\ge 1$.
If we choose $\mu,\epsilon>0$ so small such that
$d_\vp+\epsilon+\mu-s(V)<-d_\Gamma$,
then the sum converges and the right-hand side can be estimated
by $C_3 a^{\rho(s(V)+1-\mu)}$ with $C_3$ independent of $a\ge 1$.
We conclude that 
$$\lim_{a\to\infty} a^{-\rho(t(V)+1)}\sum_{g\in \Gamma}
\int_{\{x\in F|a\le a(gx)\le a_0a\} }\langle \vp(g)P f(x),\tilde h(gx)\rangle
dx=0\ ,$$ and thus $\langle f,h\rangle =0$. 

It  remains to prove    Lemma \ref{lemm3} and Lemma  \ref{claimm}.

\proof[of Lemma \ref{lemm3}]
Note that for all $g\in \Gamma$ one of the following two conditions
fails:
\begin{eqnarray*}
gg_0^{i}\cO&\in& X_+\\
g F\cap  (X\setminus X_+)&\not=&\emptyset \ .
\end{eqnarray*}
Indeed, if the first condition holds,
then $gF\cap   X_+\not=\emptyset$.
We conclude that $gF\subset    X_+$ and hence
$g F\cap  (X\setminus X_+) =\emptyset$. Further note that 
$$\overline{\{g g_0^i\cO| g\in\Gamma \:\mbox{and}\: g F\cap 
(X\setminus X_+)\not=\emptyset\}}\cap \partial X = \overline{\{g \cO|
g\in\Gamma \:\mbox{and}\: g F\cap  (X\setminus X_+)\not=\emptyset\}}\cap
\partial X\ .$$ We see that $\bar{F}\cap \partial X\subset \interi\partial X_+$
and $\overline{\{g \cO| g\in\Gamma \mbox{and} g F\cap  (X\setminus
X_+)\not=\emptyset\}}\cap \partial X\subset \partial X\setminus \partial X_+$
are disjoint. We now obtain the desired inequality from   Corollary 2.5 of
\cite{bunkeolbrich982}. \hB

\proof[of Lemma \ref{claimm}]

It is at this point where we use that $\Gamma$ is geometrically finite.
Namely, let $F_0\cup F_1\cup\dots
F_r$ be the decomposition of $F$ given in Lemma \ref{dec}.
 
Since $\Lambda_\Gamma$ and $\bar F_0$ are separated we can use  
\cite{bunkeolbrich982}, Lemma 6.2 (2), in order to get the estimate
$$\int_{\{x\in X|b\le a(x)\le a_1 b\}\cap F_0} |Pf(x)| dx\le C
b^{\rho(1-s(V))}\ ,$$
where $C$ is independent of $b\ge 1$.
This is the required estimate for the contribution of $F_0$. 

It remains to consider the contributions of the cusps, i.e of  $F_i$, $i> 0$.
Let now $[P]_\Gamma$, $P=P_i$ for one $i>0$, be a cusp of $\Gamma$. 
Then for $v\in V_{\tilde\vp}$ and $x\in X$ we have 
$$\langle Pf(x),v\rangle = \sum_{g\in \Gamma_P} \langle \vp(g)
P(\chi^{\Gamma_P} f_{|\Omega_{\Gamma_P}})(g^{-1}x), v\rangle\ .$$  Indeed, let 
$p_{x,v}\in C^\infty(\partial X,\tilde V)\otimes V_{\tilde \vp}$ denote the
integral kernel of the map $f\mapsto \langle P(f)(x), v\rangle $. Then
using the invariance properties of the kernel $p_{x,v}$ and that $f$ is
strongly supported on the limit set we get 
\begin{eqnarray}
\langle Pf(x),v\rangle &=&\langle p_{x,v},f\rangle\nonumber\\ 
&=&\sum_{g\in \Gamma_P} \langle
\chi^{\Gamma_P}f_{|\Omega_{\Gamma_P}}, \tilde \vp(g)^{-1} g^*p_{x,v}\rangle\nonumber\\ 
&=&\sum_{g\in
\Gamma_P} \langle \chi^{\Gamma_P} f_{|\Omega_{\Gamma_P}},p_{g^{-1}x,\tilde \vp(g)v}
\rangle\nonumber\\ &=&\sum_{g\in \Gamma_P} \langle \vp(g)^{-1}
P(\chi^{\Gamma_P} f_{|\Omega_{\Gamma_P}})(gx), v\rangle\nonumber\ . \end{eqnarray}

Since $\overline{\Gamma_{P} F_i}\cap \supp(\chi^{\Gamma_{P_i}}f_{|\Omega_{\Gamma_P}})=\emptyset$ we
again apply  \cite{bunkeolbrich982}, Lemma 6.2 (2),  in order to get the
estimate 
\begin{equation}\label{wunsch}|P(\chi^{\Gamma_P}f_{|\Omega_{\Gamma_P}})(gx)|\le C
a(gx)^{-\rho(s(V)+1)}\ ,\end{equation} 
where $C$ is independent of $x\in F_i$
and  $g\in \Gamma_P$. In order to estimate the sum over $g\in\Gamma_P$ we  need
the following geometric lemma. 
\begin{lem}\label{geim}
There is a constant $a_3\in \R_+$ such
that for all $x\in F$ and $g\in \Gamma_P$ we have
$a(gx)\ge a_3 \max(a(g)a(x)^{-1},a(x))$.
\end{lem}
Let us postpone the proof of the lemma und continue with the estimates.
We choose $\nu>0$ sufficiently small such that
 $d_{\Gamma_P}+d_\vp-s(V)+\nu<0$.
For those $\nu$ using  Lemma \ref{geim} and (\ref{wunsch})
we obtain  for all $x\in F_i$
\begin{eqnarray*}
|P(f)(x)|&\le&C \sum_{g\in \Gamma_P} \|\vp(g)^{-1}\| a(gx)^{-\rho(1+s(V))}\\ 
&\le &C_1 \sum_{g\in \Gamma_P} a(g)^{\rho d_\vp} a(gx)^{-\rho(1+s(V))}\\
 &\le&C_2 \sum_{g\in \Gamma_P}
a(g)^{-\rho(1+d_{\Gamma_P}+\nu)}a(x)^{\rho(-s(V)+1+2d_{\Gamma_P}+2d_\vp+2\nu)}\\
&\le& C_3 a(x)^{\rho(-s(V)+1+2d_{\Gamma_P}+2d_\vp+2\nu)} 
\end{eqnarray*} 
Since $s(V)>d_\vp+d_{\Gamma_P}+1$ we can choose
$\kappa,\nu,\mu>0$ such that
$-s(V)+1+2d_{\Gamma_P}+2d_\vp+2\nu+2+\kappa <s(V)+1-\mu$.
Then
$$|P(f)(x)|\le C_3 a(x)^{\rho(1+s(V)-\mu)} a(x)^{-\rho(2+\kappa)}\ .$$
Lemma \ref{claimm} now follows from the fact that the function $X\ni
x\mapsto a(x)^{-\rho(2+\kappa)}\in \R$ is integrable. In fact,\
$$\int_{\{x\in F_i|b\le a(x)\le a_1 b\} } |Pf(x)| dx\le C_3
b^{\rho(1+s(V)-\mu)}\int_X a(x)^{-\rho(2+\kappa)} dx
\le C_4 b^{\rho(1+s(V)-\mu)}\ ,$$
where $C_4$ is independent of $b$.
\hB

\proof[of Lemma \ref{geim}]
We consider the triangle inequality for the triangle $(\cO,gx,g\cO)$ in $X$
and obtain $\dist(\cO,gx )+\dist(g\cO,gx )\ge \dist(\cO,g\cO)$.
Since $\dist(g\cO,gx ) =\dist(\cO,x )$ we conclude that
$a(x) a(gx) \ge a(g)$. 

Recall that
$F$ is a Dirichlet domain of $\Gamma$ with respect to $g_0^i \cO$.
We conclude that $\dist(g_0^i \cO,x)\le  \dist(g_0^i \cO,gx)$
for all $g\in \Gamma_P$ and $x\in F$.
Using this and the triangle inequality for
the triangle $(\cO,x,g_0^i\cO)$ we obtain
\begin{eqnarray*}
\dist(\cO,x)& \le& \dist(g_0^i \cO,x) + \dist (\cO,g_0^i\cO)\\
&\le&\dist(g_0^i \cO,gx) + \dist (\cO,g_0^i\cO)\\
&\le&2\dist(g_0^i \cO,\cO)+\dist(\cO,gx)
\end{eqnarray*}
and therefore $a(gx)\ge a(g_0^i)^{-2}  a(x)$
for all $g\in \Gamma_P$ and $x\in F$.
The assertion of the lemma holds true if we set $a_3:=\min(1, a(g_0^i)^{-2})$
 \hB

We now have finished the proof of Lemma \ref{vam1} and therefore of
Proposition \ref{poo}. \hB 

\proof[end of the proof of Theorem \ref{main}]
By Proposition \ref{poo} we know that Theorem \ref{main} holds true under the
additional asumption that $V$ is spherical. We twist by finite-dimensional
representations of $G$ in order to conclude the general case.
Fix any parabolic subgroup $P\subset G$. Then we can write $\partial X=G/P$.
If $(\pi,V_\pi)$ is a finite-dimensional representation of $G$, then we can
form the $G$-equivariant bundle $V(\pi)=G\times_P V_\pi$ on $\partial X$.
The idea of twisting is based on the fact that there is an $G$-equivariant
isomorphism
$$C^\infty(\partial X,V\otimes V(\pi))\otimes
V_\vp\stackrel{\cong}{\rightarrow}  C^\infty(\partial X,V)\otimes V_\pi\otimes
V_\vp\ .$$ In particular, there is an isomorphism
$$j:I_{\Lambda_\Gamma}(\Gamma,V\otimes V(\pi),\vp)\stackrel{\cong}{\rightarrow}
I_{\Lambda_\Gamma}(\Gamma,V,\vp\otimes \pi)\ .$$
If $V$ is an irreducible $G$-equivariant bundle on $\partial X$, then there
exists an irreducible representation $(\pi,V_\pi)$ of $G$ and
a $G$-equivariant embedding
$$i:V\hookrightarrow \Lambda^{\frac{1-t(V)-d_\pi}{2}}\otimes V(\pi)\ ,$$
where $t(V)\in\C$ is defined such that $V\otimes \Lambda^{t(V)}=\tilde V$.
In particular, $\Ree(t(V))=s(V)$. For these facts we refer to
\cite{bunkeolbrich982}, p. 108 ( in particular,  to the formulas (33), (34)).
The embedding $i$ composed with the isomorphism $j$ gives an embedding
$$j\circ i:I_{\Lambda_\Gamma}(\Gamma,V,\vp)\hookrightarrow 
I_{\Lambda_\Gamma}(\Gamma,\Lambda^{\frac{1-t(V)-d_\pi}{2}},\vp\otimes\pi)\ .$$
We can apply Prop. \ref{poo} to the right-hand side.
Indeed,
$$s(\Lambda^{\frac{1-t(V)-d_\pi}{2}})=s(V)+d_\pi>
d_\Gamma+d_\vp+d_\pi+\max_{[P]_\Gamma }(0,d_{\Gamma_P}-d_\Gamma+1)
=d_\Gamma + d_{\vp\otimes \pi}
+\max_{[P]_\Gamma }(0,d_{\Gamma_P}-d_\Gamma+1)\ ,$$
and hence
$I_{\Lambda_\Gamma}(\Gamma,\Lambda^{\frac{1-t(V)-d_\pi}{2}},\vp\otimes\pi)=0$.
This implies $I_{\Lambda_\Gamma}(\Gamma,V,\vp)=0$. \hB

\bibliographystyle{plain}

{\small
\centerline{Ulrich Bunke and
Martin Olbrich}
\centerline{Universit\"at G\"ottingen
}\centerline{Mathematisches Institut}\centerline{Bunsenstr.
3-5}\centerline{37073
G\"ottingen}\centerline{GERMANY}\centerline{bunke@uni-math.gwdg.de,$\:\:$olbrich@uni-math.gwdg.de}
}
\end{document}